\documentclass{birkjour}
\newtheorem{theorem}{Theorem}[section]
\newtheorem{lemma}[theorem]{Lemma}
\newtheorem{corollary}[theorem]{Corollary}
\newtheorem{proposition}[theorem]{Proposition}

\theoremstyle{definition}
\newtheorem{definition}[theorem]{Definition}
\newtheorem{example}[theorem]{Example}
\usepackage{hyperref}
\theoremstyle{remark}
\newtheorem{remark}{Remark}

\begin{document}

%
%
%
%
%
%
%
%
%

\title[The g$n$s and p$n$s-Drazin inverse of sum of elements in Banach algebras]
 {The generalized and pseudo $n$-strong Drazin inverse of the sum of elements in Banach algebras}

\author{Rounak Biswas}

\address{%
 National Institute of Technology Karnataka, Surathkal, India}

\email{xyzrounak3@gmail.com}

\author{Falguni Roy}
\address{ National Institute of Technology Karnataka, Surathkal, India}
\email{falguni.roy.92@gmail.com}


\date{August 2, 2023}

\begin{abstract}
In this paper, we begin by introducing some necessary and sufficient conditions for generalized $n$-strong Drazin invertibility (g$n$s-invertibility) and pseudo $n$-strong Drazin invertibility (p$n$s-invertibility) of an element in a Banach algebra for $n\in\mathbb{N}$. Subsequently, these results are utilized to prove some additive properties of g$n$s (p$n$s)-Drazin inverse in a Banach algebra. This process produces a generalization of some recent results of 
H Chen, M Sheibani (Linear and Multilinear Algebra \textbf{70.1} (2022): 53-65) for g$n$s and p$n$s-Drazin inverse. Furthermore, we define and characterize weighted g$n$s and weighted p$n$s-Drazin inverse in a Banach algebra.
\end{abstract}

\maketitle

\section{Introduction \label{ga}}	

Throughout this paper, $\mathcal{A}$ denote a Banach algebra with the unit $1.$ The symbols $\rho(a)$, $\sigma(a)$ represents the resolvent and the spectrum of an element $a\in \mathcal{A}$, respectively. $\mathcal{A}^{qnil}$ denote the collection of all quasinilpotents elements of $\mathcal{A}$, i.e.
\begin{equation*}
\mathcal{A}^{qnil} =\{a \in \mathcal{A} : \lim_{n\to\infty} \|a^n\|^{\frac{1}{n} }=0\}=\{a \in \mathcal{A} : \sigma(a)=0\}.    
\end{equation*}

\noindent The commutant and the double commutant for $a\in \mathcal{A}$ are define by  $\text{Comm}(a)= \{x \in \mathcal{A}: xa=ax\}$ and Comm$^2(a)=\{x\in \mathcal{A} : xy =yx ,\hspace{0.2cm} \forall \hspace{0.2cm} y\in \text{Comm}(a)\}$. Set of all invertible and nilpotent elements of $\mathcal{A}$ will be denoted by $\mathcal{A}^{inv}$ and $\mathcal{A}^{nil}$, respectively.
 In 1958, Drazin \cite{drazin1958pseudo} introduced a new kind of generalized inverse different from the well-known Moore-Penrose generalized inverse, known as the Drazin inverse. An element  $a\in\mathcal{A}$ is Drazin invertible if there exists $b\in\mathcal{A}$ such that $$ab=ba,\text{ }b=bab\text{ and }a-a^2b\in\mathcal{A}^{nil}.$$ 
  The generalized Drazin invertibility (g-Drazin invertibility) of an element $a \in \mathcal{A}$ is defined in \cite{koliha1996generalized} by replacing $a-a^2b\in\mathcal{A}^{nil}$ with $a-a^2b\in\mathcal{A}^{qnil}$. Hence, the g-Drazin inverse of $a$ is the element $b \in \mathcal{A}(\text{written }a^d)$ which satisfies $$ab=ba\text{, } b=ab^2\text{, }a-a^2b \in A^{qnil}.$$ The set of all g-Drazin invertible elements in $\mathcal{A}$ will be denoted by $\mathcal{A}^d$ and $a^{\pi}=1-aa^d$, known as the spectral idempotent of $a.$  
  The theory of Drazin inverses has grown significantly over the past few decades, and it
  finds applications in many areas, including Statistics, Numerical analysis, Differential
equations, Population models, Cryptography, and Control theory,
etc.

    \noindent Wang \cite{wang2017class} introduced a new subclass of Drazin invertible elements, known as strong Drazin invertible elements equivalent to the notion of strongly nil-clean elements presented in \cite{diesl2013nil}. An element $a\in\mathcal{A}$ is strong Drazin invertible, if there exists $b\in\mathcal{A}$ such that $$ab=ba,\text{ }bab=b\text{ and }a-ab\in\mathcal{A}^{nil}.$$ Motivated by g-Drazin inverse, a modification of strong Drazin inverse was presented in \cite{mosic2016reverse}, as  generalized strong Drazin inverse(gs-Drazin inverse). The gs-Drazin inverse of $a\in\mathcal{A}$ is the element $x\in \mathcal{A}$ which satisfies $$xax=x,\text{ }ax=xa,\text{, }a-ax\in \mathcal{A}^{qnil}.$$ Like strong Drazin invertible elements, the set of all gs-Drazin invertible elements in $\mathcal{A}$ form a subclass of g-Drazin invertible elements in $\mathcal{A}$.
  \\ Given $n \in \mathbb{N}$, Dijana Mosić \cite{mosic2021generalized} introduced the notion of generalized $n$-strong invertibility in a Ring. In a Banach algebra, an element $a \in \mathcal{A}$ is  generalized $n$-strongly Drazin invertible (g$n$s-Drazin invertible) if there exists an element $x \in \mathcal{A}$ such that $$x \in \text{Comm}^2(a),\text{ } xax=x\text{ } \text{and}\text{ } a^n-ax \in \mathcal{A}^{qnil}.$$ Then the element $x$ is the generalized $n$-strong Drazin inverse of $a$ and denoted by $a^{nsd}$. If exists, $a^{nsd}$ is unique and coincide with $a^{d}$. Collection of all g$n$s-Drazin invertible elements in $\mathcal{A}$ is denoted by $\mathcal{A}^{nsd}$. In particular, for $n=1$, g$n$s-Drazin inverse is the gs-Drazin inverse, and for $n=2$, $2$sd-Drazin inverse is known as generalized Hirano inverse \cite{chen2019generalized}. The collection of generalized Hirano invertible elements in $\mathcal{A}$ will be denoted by $\mathcal{A}^H$, and if $a\in \mathcal{A}$ is generalized Hirano invertible, then the corresponding generalized Hirano inverse is denoted by $a^H.$
  
\noindent Zhou Wang and Jianlong Chen \cite{wang2012pseudo} presented an intermediary inverse between Drazin and g-Drazin inverse, known as pseudo Drazin inverse. This particular inverse is obtained by replacing the quasinilpotency condition in g-Drazin inverse definition with Jacobson radical.
The Jacobson radical $\mathcal{J}(\mathcal{A})$ in a Banach algebra $\mathcal{A}$ is defined as 
\begin{equation*}
    \mathcal{J}(\mathcal{A})=\{a\in \mathcal{A} : 1+ay \in \mathcal{A}^{inv}, \text{ for all } y \in\mathcal{A}\},
\end{equation*}
and $$\sqrt{\mathcal{J}(\mathcal{A})}=\{a\in \mathcal{A}:a^n \in \mathcal{J}(\mathcal{A})\text{ for some }n \in \mathbb{N}\}.$$ An element $a\in \mathcal{A}$ is pseudo Drazin invertible(p-Drazin invertible) if there exists $b\in\mathcal{A}$ such that 
    \begin{equation*}
        b\in\text{Comm}^2(a),\text{ }ab^2=b\text{ and }a^k(1-ab)\in \mathcal{J}(\mathcal{A})\text{ for some }\hspace{.2cm}k\geq1.
    \end{equation*}
    Here $b$ is the p-Drazin inverse of $a$, denoted by $a^{pD},$  and
    the collection of all p-Drazin invertible elements of $\mathcal{A}$ will be denoted by $\mathcal{A}^{pD}.$
  \\ For $n\in\mathbb{N}$, Mosic \cite{mosic2021generalized} presented the concept of pseudo $n$-strong Drazin inverse(p$n$s-Drazin inverse). The p$n$s-Drazin inverse of $a\in \mathcal{A}$ is the element $x\in\mathcal{A}(\text{written }a^{pnsD})$ such that 
    \begin{equation*}
        x\in\text{Comm}^2(a),\text{ }xax=x\text{ and } (a^n-ax)^k\in\mathcal{J}(\mathcal{A}) \text{ for some } k\in\mathbb{N}.
    \end{equation*}
    The set of all p$n$s-Drazin invertible elements of $\mathcal{A}$ will be denoted by $\mathcal{A}^{pnsD}$. Mosic \cite{mosic2021generalized} also defined weighted generalized strong Drazin inverse(wgs-Drazin inverse) and weighted pseudo strongly Drazin inverse(wps-Drazin inverses) for some non-zero weight in Banach algebra. 
  \\ One of the most interesting problems in g-Drazin inverse is the study of g-Drazin invertibility for the sum of two g-Drazin invertible elements. Under non-identical circumstances, several studies \cite{castro2010expressions,cvetkovic2009generalized,deng2010new,chen2022g} have proposed their work on this problem. Let $a,b\in\mathcal{A}^d.$ If $a^{\pi}b=b,\text{ }ab^{\pi}=a,\text{ }b^{\pi}aba^{\pi}=0$, then it is proved in \cite{gonzalez2004new} that $a+b$ is g-Drazin invertible. Cvetković-Ilić et al. \cite{cvetkovic2006additive},  proved that if $a=ab^{\pi},\text{ }b^{\pi}ba^{\pi}=b^{\pi}b\text{ and }b^{\pi}a^{\pi}ba=b^{\pi}a^{\pi}ab$ then $a+b\in\mathcal{A}^d$.
  
  \noindent Recently, some additive properties for gs-Drazin and n-strong Drazin inverse are studied in \cite{chen2020gs} and \cite{zou2019n}, respectively.  Motivated by these works, we provide some equivalent conditions for the g$n$s and p$n$s-Drazin invertibility of an element in a Banach algebra  with some new additive results, which extend the works of \cite{chen2020gs}, \cite{zou2019n} and \cite{chen2022g}.  
  \\  In Section \ref{sec1}, we investigate the g$n$s-Drazin invertibility of an element in a Banach algebra and derive some necessary and sufficient conditions for g$n$s-Drazin invertibility. Subsequently, based on these results, we established several additive properties for g$n$s-Drazin inverse. In the next section, we look over the p$n$s-Drazin invertibility of an element in Banach algebra and provide some additive results for the p$n$s-Drazin inverse. Finally, in Section \ref{sec4}, we define and characterize the weighted g$n$s and p$n$s-Drazin invertibility in a Banach algebra. 

\section{Additive results for g$n$s-Drazin inverse. \label{sec1}}
\indent In this section, we give some necessary and sufficient conditions for the g$n$s-Drazin invertibility of an element in a Banach algebra. Then, we prove a few additive results for the g$n$s-Drazin inverse.  
\\We begin with some results on quasinilpotent elements of a Banach algebra.
\begin{lemma}\label{l1} \cite{cvetkovic2009generalized}
Let $a,b \in \mathcal{A}$ and $ab=\lambda ba$ for some scalar $\lambda$ and $ab\neq 0$.
\begin{enumerate}
    \item If $a,b\in \mathcal{A}^{qnil}$ then $a+b \in \mathcal{A}^{qnil}$;
    \item If $a$ or $b$ $\in \mathcal{A}^{qnil}$ then $ab\in \mathcal{A}^{qnil}$.
\end{enumerate}
\end{lemma}
\begin{lemma}\label{l2} \cite{gonzalez2004new}
If $a, b \in \mathcal{A}^{qnil}$ and $ab=0$, then $a+b \in \mathcal{A}^{qnil}$.

\end{lemma}
Now, we establish a key proposition for our main theorem. 

\begin{proposition}\label{p1}
Let $n\in\mathbb{N}$, and $a \in\mathcal{A}$. Then $a^n-a^{2n}\in \mathcal{A}^{qnil}$ if and only if $a-a^{n+1} \in \mathcal{A}^{qnil}$.
\begin{proof}
Assume that $a^n-a^{2n} \in \mathcal{A}^{qnil}$. Since we have\begin{align*}
(a-a^{n+1})^n&=(a-a^{n+1})a^{n-1}(1-a^n)^{n-1}\\&=(a^n-a^{2n})(1-a^n)^{n-1}.
\end{align*}Then $(a-a^{n+1})^n\in \mathcal{A}^{qnil}$ by Lemma \ref{l1}, therefore $(a-a^{n+1})\in\mathcal{A}^{qnil}.$
\\ \indent Conversely, if $a-a^{n+1}\in \mathcal{A}^{qnil}$, then $a^n-a^{2n}=a^{n-1}(a-a^{n+1})\in \mathcal{A}^{qnil}$, using Lemma \ref{l1}.  
\end{proof}
\end{proposition}
We come now to the presentation of the main result of this section.
\begin{theorem}\label{t1}
Let $a\in\mathcal{A}$ and $n\in \mathbb{N}$. Then $a$ is g$n$s-Drazin invertible if and only if $a^n-a^{2n}\in\mathcal{A}^{qnil}.$
    \begin{proof}
    Since $a\in \mathcal{A}^{nsd}$, then we have $a^n-aa^{nsd} \in \mathcal{A}^{qnil}$.
    Now
    \begin{align*}
       a^n-a^{2n}&=a^n-aa^{nsd}+aa^{nsd}-a^{2n}\\&=(a^n-aa^{nsd})+(aa^{nsd})^2-a^{2n}\\&=(a^n-aa^{nsd})+(aa^{nsd}-a^n)(aa^{nsd}+a^n). 
    \end{align*}
    Hence by Lemma \ref{l1}, we obtain $a^n-a^{2n} \in  \mathcal{A}^{qnil}$.
    \\Conversely, if $a^n-a^{2n} \in  \mathcal{A}^{qnil}$, then by Proposition \ref{p1}, $a-a^{n+1}\in  \mathcal{A}^{qnil}$, which implies  $\sigma(a-a^{n+1})=\{0\}$. Now consider the polynomial $p(x)=x-x^{n+1}$ over $\mathbb{C}$. Then using Spectral Mapping Theorem \cite[Theorem 4.10]{conway2019course} we have $$\{p(\lambda): \lambda \in \sigma(a)\} = \sigma(p(a))=\sigma(a-a^{n+1})=\{0\},$$ which implies $\sigma(a)\subset \{\lambda \in \mathbb{C}:\lambda^{n}=1\hspace{0.2cm}\text{or}\hspace{0.2cm}  \lambda=0\}$. Thus $0$ can only be an isolated spectral point of $a$. Hence $a$ is g-Drazin invertible with the g-Drazin inverse $a^d$. Since $a-a^2a^d\in  \mathcal{A}^{qnil}$ then $$(a^n-1+aa^d)((a^d)^n-1+aa^d)=1-(a^n-a^{n+1}a^d)\in \mathcal{A}^{inv}.$$ Therefore in particular $a^n-1+aa^d$ is invertible. Moreover $$(a^n-aa^d)(a^n-1+aa^d)=a^{2n}-a^n \in  \mathcal{A}^{qnil},$$ then $(a^n-aa^d)=(a^n-1+aa^d)^{-1}(a^{2n}-a^n) \in  \mathcal{A}^{qnil}$ by Lemma \ref{l1}. Hence $a\in  \mathcal{A}$ is g$n$s-Drazin invertible.
    \end{proof}

\end{theorem}
\begin{corollary}\label{c1}
Let $a\in \mathcal{A}$ and $n\in\mathbb{N}$. Then $a$ is gns-Drazin invertible if and only if $a-a^{n+1}\in\mathcal{A}^{qnil}$.
\begin{proof}
It follows from Theorem \ref{t1} and Proposition \ref{p1}.
\end{proof}

\end{corollary}
\begin{remark}
    Chen and Calci \cite{chen2020gs}, proved $a\in\mathcal{A}$ is gs-Drazin invertible if and only if $a\in \mathcal{A}^d$ and $a-a^2\in\mathcal{A}^{qnil}$ but by the Corollary \ref{c1}, only second condition is enough.
\end{remark}
\begin{corollary}\label{c2}
Let $a\in\mathcal{A}$ and $n \in \mathbb{N}$.
\begin{enumerate}
    \item \label{new1} If $a\in \mathcal{A}^{sd}$, then $a\in \mathcal{A}^{nsd}$ and $a^{nsd}=a^{sd}=a^d$;
    \item If $a \in \mathcal{A}^{H}$, then $a \in \mathcal{A}^{2nsd}$ and $a^{2nsd}=a^{H}=a^{d}$;
    \item If $a\in \mathcal{A}^{nsd}$, then $a \in \mathcal{A}^{2nsd}$ and $a^{2nsd}=a^{nsd}=a^{d}$.
\end{enumerate}
\begin{proof}
\begin{enumerate}
    \item Since $a\in \mathcal{A}^{sd}$, therefore by Corollary \ref{c1}, $a-a^2 \in  \mathcal{A}^{qnil}$. Now $$a-a^{n+1}=a(1-a) \displaystyle(\sum_{i=0}^{n-1}a^i)=(a-a^2)(\sum_{i=0}^{n-1}a^i)\in \mathcal{A}^{qnil}.$$ Hence $a\in \mathcal{A}^{nsd}$ and by uniqueness of gs-Drazin inverse $a^{nsd}=a^{sd}=a^{d}$. 
    \item If $a\in\mathcal{A}^H$, then by the Corollary \ref{c1}, $a-a^3\in \mathcal{A}^{qnil}$. Now $a-a^{2n+1}=(a-a^3)(1+a^2+a^4\dots+a^{2n-2})\in \mathcal{A}^{qnil}$ , therefore the required result follows from Corollary \ref{c1}.
    \item This follows in a similar manner as \ref{new1}.
\end{enumerate}
\end{proof}
\end{corollary}
Koliha\cite{koliha1996generalized} proved the g-Drazin invertibility for the product of two commutative g-Drazin invertible elements in a Banach algebra. Using Theorem \ref{t1}, we will prove the g$n$s-Drazin invertibility for the product of two g$n$s-Drazin invertible elements in a Banach algebra under a condition weaker than commutativity.
\begin{lemma}\label{c3}
Let $n\in\mathbb{N}$ and  $a,b\in\mathcal{A}^{nsd}$, if  $a^2b=aba$. Then $ab\in\mathcal{A}^{nsd}$.
\begin{proof} Since $a^2b=aba$, then $(ab)^m=a^mb^m$ for any $m\in\mathbb{N}.$ 
From Theorem \ref{t1}, we have $$a-a^{n+1},b-b^{n+1}\in \mathcal{A}^{qnil}.$$ Moreover
\begin{equation}\label{eqn1}
    (ab)-(ab)^{n+1}=a^{n+1}(b-b^{n+1})+(a-a^{n+1})b,
\end{equation}
where $$(a-a^{n+1})ba^{n+1}(b-b^{n+1})=a^{n+1}(b-b^{n+1})(a-a^{n+1})b.$$
 Therefore in light of  Lemma \ref{l1}, $(ab)-(ab)^{n+1}\in\mathcal{A}^{qnil}$. Hence by Corollary \ref{c1}, $ab\in \mathcal{A}^{nsd}$. 
\end{proof}
\end{lemma}
\begin{corollary}
    Let $n\in\mathbb{N}$. If $a,b\in\mathcal{A}^{nsd}$ and $ab=ba$. Then $ab\in\mathcal{A}^{nsd}.$
    \begin{proof}
        Follows from Lemma \ref{c3}.
    \end{proof}
\end{corollary}
In the following theorems, we investigate the g$n$s-Drazin invertibility of the sum of two g$n$s-Drazin invertible elements in Banach algebra under different conditions. 
\begin{theorem}\label{t2}
Let $n\in\mathbb{N}$ and $a, b\in \mathcal{A}^{nsd}$. If $ab=0$, then $a+b$ is gns-Drazin invertible.
\begin{proof}
If $a,b\in \mathcal{A}^{nsd}$, then by Corollary \ref{c1}, we have $a-a^{n+1}$ and $b-b^{n+1}$ $\in\mathcal{A}^{qnil}$. Since $$(a+b)-(a+b)^{n+1}=(a-a^{n+1})+(b-b^{n+1})+(\sum_{i=1}^{n}b^ia^{n+1-i}),$$ therefore using Lemma \ref{l2},  we obtain $(a+b)-(a+b)^{n+1}\in \mathcal{A}^{qnil}$. Thus $a+b\in \mathcal{A}^{nsd}$ by Corollary \ref{c1}. 
\end{proof}
\end{theorem} 
\begin{theorem}\label{theo}
    Let $n\in\mathbb{N}$ and $a,b\in\mathcal{A}^{nsd}$ such that $ab= ba$, if $ab\in\mathcal{A}^{qnil}$ then $a+b\in\mathcal{A}^{nsd}$ and if $(a+b)\in \mathcal{A}^{nsd}$, then $ab\left(\displaystyle\sum_{k=1}^n\binom{n+1}{k}a^{k-1}b^{n-k}\right)\in \mathcal{A}^{qnil}.$
    \begin{proof}
        By Corollary \ref{c1}, $a+b\in \mathcal{A}^{nsd}$ if and only if $(a+b)-(a+b)^{n+1}\in \mathcal{A}^{qnil}.$ Moreover $$(a+b)-(a+b)^{n+1}=(a-a^{n+1})+(b-b^{n+1})-ab\left(\sum_{k=1}^n\binom{n+1}{k}a^{k-1}b^{n-k}\right).$$ Since both $a,b$ are g$n$s-Drazin invertible, then $a-a^{n+1},b-b^{n+1}\in\mathcal{A}^{qnil}$. Therefore by Lemma \ref{l1}, $(a+b)$ is g$n$s-Drazin invertible if and only if $\sum_{k=1}^n\binom{n+1}{k}a^kb^{n+1-k}\in\mathcal{A}^{qnil}$. Now $$ab\in\mathcal{A}^{qnil}\implies ab\left(\sum_{k=1}^n\binom{n+1}{k}a^{k-1}b^{n-k}\right)\in\mathcal{A}^{qnil}.$$
    \end{proof}
\end{theorem}
\begin{remark}
    It is important to note that in Theorem \ref{theo}, $a+b\in\mathcal{A}^{nsd}$ does not necessarily implies $ab\in\mathcal{A}^{qnil}$, which can be seen in the next example.
\end{remark}
\begin{example}
    $$A=\begin{bmatrix}
        1 & 0\\
        0 & 1\\
    \end{bmatrix}\hspace{0.2cm}\text{and}\hspace{0.2cm}B=\begin{bmatrix}
        -1 & 0 \\
        0 & -1
    \end{bmatrix}\in M_2(\mathbb{R})^{2sd}.$$ Then $AB=BA$ and $A+B\in M_2(\mathbb{R})^{2sd}$, but $AB\notin M_2(\mathbb{R})^{qnil}$, where as $3AB(A+B)\in M_2(\mathbb{R})^{qnil}.$
\end{example}
If $p\in \mathcal{A}$ is an idempotent element and $y\in \mathcal{A}$, then for the total system of idempotent $\mathcal{P}=(p,1-p)$, $y$ has a unique representation as discussed in  \cite{cvetkovic2006additive},  given by
\begin{equation}\label{eq2}
y=\begin{bmatrix}
pyp & py(1-p) \\
(1-p)yp & (1-p)y(1-p)\\
\end{bmatrix}_p.
\end{equation}
 In particular for $a\in  \mathcal{A}^{d}$, the unique block representation \ref{eq2}, of $a$ corresponding to the idempotent $p=aa^d$ is $\begin{bmatrix}
    a_1 &0\\
    0 & a_2\\
\end{bmatrix}_p,$ where $a_1\in p\mathcal{A}p^{inv}$ and $a_2\in(1-p)\mathcal{A}(1-p)^{qnil}.$ From now onward the Banach algebra $p\mathcal{A}p$ and $(1-p)\mathcal{A}(1-p)$ will be denoted by $\mathcal{A}_1$ and $\mathcal{A}_2$, respectively. 
\begin{lemma}\label{l3}
Let $x,y \in \mathcal{A}$ such that
$x=\begin{bmatrix}
a & 0\\
c & b\\
\end{bmatrix}_p$
and $y=\begin{bmatrix}
b & c\\ 0 & a\\
\end{bmatrix}_{1-p}$

be the representation of $x$ and $y$ relative to a total system of idempotent $(p,1-p)$ and $(1-p,p)$, respectively. Then for $n\in \mathbb{N}$
\begin{enumerate}
    \item \label{111}if $a\in \mathcal{A}_1^{nsd}$ and $b\in \mathcal{A}_2^{nsd}$, then $x,y$ are g$n$s-Drazin invertible;
    \item if $x\in \mathcal{A}^{nsd}$ and $b \in {\mathcal{A}_2^{nsd}}$ then $a$ is gns-Drazin invertible in $\mathcal{A}_1$.
    \end{enumerate}
    \begin{proof}
    \begin{enumerate}
        \item Since, $a\in \mathcal{A}_1$ and $b\in \mathcal{A}_2$ are g$n$s-Drazin invertible therefore $\sigma_1(a-a^{n+1})=0$ and $\sigma_2(b-b^{n+1})=0$, where $\displaystyle\sigma_1(a-a^{n+1})$ and $\sigma_2(b-b^{n+1})$ denotes their spectrum in $\mathcal{A}_1$ and $\mathcal{A}_2$, respectively. Moreover
        $$x-x^{n+1}=\begin{bmatrix}
        a & 0 \\c & b\\
        \end{bmatrix}_p-\begin{bmatrix}
        a^{n+1} & 0 \\ \alpha_{n+1} & b^{n+1}\\
        \end{bmatrix}_{p}=\begin{bmatrix}
        a-a^{n+1} & 0\\
        c-\alpha_{n+1} & b-b^{n+1}
        \end{bmatrix}_p,$$ where $\hspace{0.2cm}\alpha_{n+1}=\displaystyle{\sum_{k=0}^{n}b^kca^{n-k}}.$
        Therefore $\sigma(x-x^{n+1}) \subseteq \sigma_1(a-a^{n+1})\cup \sigma_2(b-b^{n+1})$, hence $x-x^{n+1}\in\mathcal{A}^{qnil}$. This gives $x$ is gns-Drazin invertible by Corollary \ref{c1}. Similarly,  it follows for $y$ also.
        \item Similar to \ref{111}, this result follows from the fact that $\sigma(a)\subset \sigma(x) \cup \sigma(b) $.
    \end{enumerate}
    \end{proof}
    \end{lemma}
    Next, we analyse the g$n$s-Drazin invertibility for the sum of two particular types of g$n$s-Drazin invertible elements in a Banach algebra under a condition weaker than commutativity. This also generalizes some results from \cite{cvetkovic2006additive}.
    
    \begin{lemma}\label{l4}
 Let $a\in \mathcal{A}^{qnil}$ and $b\in \mathcal{A}$ is g$n$s-Drazin invertible for some $n\in\mathbb{N}$. If $a=ab^{\pi}$, $b^{\pi}ba=b^{\pi}ab$, then $a+b$ is g$n$s-Drazin invertible.
\begin{proof}
Since $b\in \mathcal{A}^{nsd}$, therefore $b\in \mathcal{A}^d$. For $p=1-b^{\pi}$, we have the representation $$b=\begin{bmatrix}
b_1 & 0\\ 0 & b_2\\
\end{bmatrix}_p \text{ and }a=\begin{bmatrix}
a_1 & a_2\\a_3 & a_4\\
\end{bmatrix}_p,$$ where $b_1 \in \mathcal{A}_1^{inv}$ is invertible and $b_2 \in \mathcal{A}_2^{qnil}$. Use of this representation on the condition $a=ab^{\pi}$ gives $$ \begin{bmatrix}
a_1 & a_2\\a_3 & a_4\\
\end{bmatrix}_p=\begin{bmatrix}
a_1 & a_2\\ a_3 & a_4\\
\end{bmatrix}_p
\begin{bmatrix}
0 & 0\\0 & 1\\
\end{bmatrix}_p\implies a_1=a_3=0.$$ Therefore we get $a=\begin{bmatrix}
0 & a_2\\ 0 & a_4\\
\end{bmatrix}_p$
 and $a+b= \begin{bmatrix}
b_1 & a_2 \\ 0 & a_4 +b_2\\
\end{bmatrix}_p$. Again $a=ab^{\pi}$ gives $a_4=b^{\pi}ab^{\pi}=b^{\pi}a$  and $a_4^n=b^{\pi}a^n$ thus $a_4 \in \mathcal{A}^{qnil} $, as $a\in \mathcal{A}^{qnil}$. Similarly, using block matrix representations of $b$ and $a$ on the other condition $b^{\pi}ab=b^{\pi}ba$ we obtain $a_4b_2=b_2a_4$. Therefore by Lemma \ref{l1}, $a_4+b_2\in \mathcal{A}^{qnil}$. Also, $b\in \mathcal{A}^{nsd} $ implies $b_1=(1-b^{\pi})b\in \mathcal{A}^{nsd}$. Hence, by Lemma \ref{l3}, $a+b \in \mathcal{A}^{nsd}$.
\end{proof}
\end{lemma}
Now we are ready to prove the following.
\begin{theorem}\label{t3}
Let $n\in\mathbb{N}$ and $a,b \in \mathcal{A}^{nsd}$. If $a=ab^{\pi}$, $b^{\pi}ba^{\pi}=b^{\pi}b$ and $b^{\pi}a^{\pi}ba=b^{\pi}a^{\pi}ab$ then $a+b\in \mathcal{A}^{nsd}$.
\begin{proof}
Let $p=1-b^{\pi}$. Then $b=\begin{bmatrix}
b_1 & 0\\ 0 & b_2\\
\end{bmatrix}_p$, where $b_1 \in \mathcal{A}_1^{inv}$ and $b_2 \in \mathcal{A}_2^{qnil}$ and let $a=\begin{bmatrix}
a_1 & a_2 \\ a_3 & a_4\\
\end{bmatrix}_p$. Using the representation of $a$ and $b$ on the condition $a=ab^\pi,$ we obtain $a_1=a_3=0$, therefore $a=\begin{bmatrix}
0 & a_2\\ 0 & a_4\\
\end{bmatrix}_p.$  Moreover $$b^{\pi}ba^{\pi}=b^{\pi}b\implies b_2a_4^{\pi}=b_2\text{ and }b^{\pi}a^{\pi}ba=b^{\pi}a^{\pi}ab \implies a_4^{\pi}b_2a_4=a_4^{\pi}a_4b_2.$$
Since  
$ b_2\in\mathcal{A}_2^{qnil}$ and $a_4=b^{\pi}ab^{\pi}=b^{\pi}a$ is gns-Drazin invertible satisfying $b_2=b_2a_4^{\pi}$,  $a_4^{\pi}a_4b_2=a_4^{\pi}b_2a_4$, therefore by Lemma \ref{l4}, $a_2+b_4\in (1-p)A(1-p)^{nsd}$. Moreover $a+b= \begin{bmatrix}
b_1 & a_2\\ 0 & b_2+a_4\\
\end{bmatrix}_p$ where $b_1=(1-b^{\pi})b\in\mathcal{A}_1^{nsd}$,
hence $a+b\in \mathcal{A}^{nsd}$ in light of Lemma \ref{l3}.
\end{proof}
\end{theorem}
We can obtain the next corollary as a special case of the preceding theorem. 
\begin{corollary}
If $a,b \in \mathcal{A}^{nsd}$ for some $n\in\mathbb{N}$ and 
$ab=ba$,$a=ab^{\pi}$, $b^{\pi}=ba^{\pi}=b^{\pi}b$, then $a+b$ is g$n$s-Drazin invertible.

\end{corollary}
  Next, we investigate the g$n$s-Drazin invertibility for the sum of two g$n$s-Drazin invertible elements in a Banach algebra under the same conditions as used in \cite{gonzalez2004new}.
\begin{lemma} \label{l5}
Let $n\in\mathbb{N}$, $b\in \mathcal{A}^{nsd}$ and $a\in \mathcal{A}^{qnil}.$ If $ab^{\pi}=a$ and $b^{\pi}ab=0$, then $a+b\in \mathcal{A}^{nsd}$.
\begin{proof}
If $b\in \mathcal{A}^{qnil}$ then $b^{\pi}=1$ therefore $ab=0$, then by Theorem \ref{t2}, $(a+b)$ is g$n$s-Drazin invertible. So assume that $b\notin \mathcal{A}^{qnil}$ then for $p=1-b^{\pi}$ we have $$b=\begin{bmatrix}
b_1 & 0\\ 0 & b_2\\
\end{bmatrix}_p\text{ and }a=\begin{bmatrix}
a_4 & a_1\\a_3 & a_2\\
\end{bmatrix}_p,$$ where $b_1 \in \mathcal{A}_1^{inv}$ and $b_2 \in \mathcal{A}_2^{qnil}$. Now using this representation of $a$ and $b$ on $a(1-b^{\pi})=0$ we obtain $a_4=a_3=0$, therefore $a=\begin{bmatrix}
0 & a_1\\0 & a_2\\
\end{bmatrix}_p$. 
Again from $b^{\pi}ab=0$, it follows that $ a_2b_2=0$. Since $a_2=b^{\pi}a\in A^{qnil}$ and  $b_2=b^{\pi}b\in A^{qnil}$, then by Lemma \ref{l2}, $a_2+b_2\in \mathcal{A}^{qnil}.$ Moreover
$b_1=(1-b^{\pi})b\in \mathcal{A}^{nsd}$, therefore using Lemma \ref{l3}, $a+b=\begin{bmatrix}
    b_1 & a_1\\
    0 & b_2+a_2\\
\end{bmatrix}_p\in \mathcal{A}^{nsd}$.
\end{proof}
\end{lemma}
We are now prepared to demonstrate the following result.
\begin{theorem}\label{t4}
Let $a,b\in\mathcal{A}^{nsd}$ for some $n\in\mathbb{N}$. If $a=ab^{\pi}, b=ba^{\pi}$ and $b^{\pi}aba^{\pi}=0$ then $a+b\in \mathcal{A}^{nsd}$.
\begin{proof}
Let $p=1-a^{\pi}$ then $a=\begin{bmatrix}
a_1& 0\\ 
0 & a_2\\
\end{bmatrix}_p$, where $a_1\in (\mathcal{A}_1)^{inv}$ and $a_2\in \mathcal{A}_2^{qnil}$ and let $b=\begin{bmatrix}
b_1 & b_2\\
b_3 & b_4\\
\end{bmatrix}_p$. Now $b=ba^{\pi}$ implies $b_1=b_3=0$, therefore $b=\begin{bmatrix}
0 & b_2\\
0 & b_4\\
\end{bmatrix}_p$. Since $b\in \mathcal{A}^{nsd}$ then by Lemma \ref{l3}, $b_4\in\mathcal{A}_2^{nsd}$ and $b^{\pi}=\begin{bmatrix}
1 & -b_2b_4^d\\
0 & b_4^{\pi}
\end{bmatrix}_p$. Furthermore using this representation of $a$ and $b^\pi$ on  $a=ab^{\pi}$ we obtain 
 $a_2=a_2b_4^{\pi}$.
Moreover $$b^{\pi}aba^{\pi}=0 \implies b_4^{\pi}a_2b_4=0.$$ Then using Lemma \ref{l5}, $a_2+b_4 \in \mathcal{A}_2^{nsd}$. Since $a+b=\begin{bmatrix}
a_1 & 0\\ b_2 & a_2+b_4\\
\end{bmatrix}_p$ where $a_1=a(1-a^{\pi})$ and $a_2+b_4$ both are gns-Drazin invertible in their respective Banach algebra therefore $a+b\in \mathcal{A}^{nsd}$ by Lemma \ref{l3}.
\end{proof}
\end{theorem}
\begin{theorem}\label{Theor}
    For $n\in \mathbb{N}$, if $a,b\in\mathcal{A}^{nsd}$ such that $aba^\pi=0$ and $(1-a^\pi)b\in\mathcal{A}^{nsd}$. Then $a+b\in \mathcal{A}^{nsd}$ if and only if $aa^{nsd}(a+b)\in\mathcal{A}^{nsd}.$
    \begin{proof}
        Here $a,b$ are g$n$s-Drazin invertible. Therefore for the idempotent $p=aa^d$, $a$ and $b$ has the representation $$a=\begin{bmatrix}
            a_1 & 0\\
            0 & a_2\\
        \end{bmatrix}_p \text{ and }b=\begin{bmatrix}
                b_1 & b_4\\
                b_2 & b_3\\
            \end{bmatrix}_p ,$$ where  and $a_1\in (\mathcal{A}_1)^{inv}$, $a_2\in\mathcal{A}_2^{qnil}.$ Moreover \begin{equation*}
            aba^\pi=0\implies b_4=0 \text{ and }a_2b_3=0.
            \end{equation*} Hence $b=\begin{bmatrix}
                b_1 & 0\\
                b_2 & b_3\\
            \end{bmatrix}_p$, where $b_3\in \mathcal{A}_2^{nsd}$ since $(1-a^\pi)b\in\mathcal{A}^{nsd}$. As $a_2\in\mathcal{A}_2^{qnil}$, $b_3\in \mathcal{A}_2^{nsd}$ and $a_2b_3=0$, then   $a_2+b_3\in\mathcal{A}_2^{nsd}$ by Theorem \ref{t2}. Hence by Lemma \ref{l3}, $a+b=\begin{bmatrix}
            a_1+b_1 & 0\\
            b_2 & a_2+b_3\\
        \end{bmatrix}_p\in \mathcal{A}^{nsd}$  if and only if $a_1+b_1\in (\mathcal{A}_1)^{nsd}$, i.e. if and only if $aa^{nsd}(a+b)\in\mathcal{A}^{nsd}.$ 
    \end{proof}
\end{theorem}
Now we establish the quasinilpotency for the sum of two quasinilpotent elements in Banach algebra under a condition weaker than commutativity.
\begin{lemma}\label{mylemma}
    If $a,b\in \mathcal{A}^{qnil}$, $ab^2=bab$ and  $a^2b=aba$, then $a+b\in \mathcal{A}^{qnil}.$
    \begin{proof}
        Since $a,b\in\mathcal{A}^{qnil},$ therefore for $0<\epsilon< \text{max}\{\|a\|,\|b\|\}$ there exists $n_0\in \mathbb{N}$ such that $$\|a^n\|<\epsilon^n \hspace{0.2cm}\text{and}\hspace{0.2cm}\|b^n\|<\epsilon^n,\text{ whenever }n>n_0.$$ 
        Now for $n>n_0$ each monomial in the expansion of $(a+b)^{3n}$ is either of the form $a^{k_1}b^{k_2}a^{k_1}$ or $b^{l_1}a^{l_2}b^{l_3}$ where $l_1+l_2+l_3=k_1+k_2+k_3=3n.$ Then $$\|(a+b)^n\|<(2\epsilon)^n\left(\text{max}\left\{\frac{\|a\|}{\epsilon},\frac{\|b\|}{\epsilon}\right\}\right)^{2n_0}.$$ Hence $a+b\in \mathcal{A}^{qnil}.$
    \end{proof}
\end{lemma}

\begin{theorem}\label{2022}
    Let $a,b\in\mathcal{A}^{nsd}$, for some $n\in\mathbb{N}$, such that $a^2b=aba$, $ab^2=bab$ and $aa^{nsd}b\in\mathcal{A}^{nsd}.$ Then  the followings are equivalent:
    \begin{enumerate}
        \item  $a+b\in \mathcal{A}^{nsd}$;\label{11}
        \item \label{it1} $aa^{nsd}(a+b)bb^{nsd}\in \mathcal{A}^{nsd}$;
        \item\label{it2} $aa^{nsd}(a+b)\in \mathcal{A}^{nsd};$
        \item\label{it3} $(a+b)bb^{nsd}\in \mathcal{A}^{nsd}.$\label{12}
    \end{enumerate}
   
    \begin{proof}
    Let $p=aa^d$, then we have $$a=\begin{bmatrix}
        a_1 & 0 \\
        0 & a_2\\
    \end{bmatrix}_p\text{ and }\begin{bmatrix}
        b_1 & b_2\\
        b_3 & b_4\\
    \end{bmatrix}_p,$$ where  $a_1\in(\mathcal{A}_1)^{nsd}$, $a_2\in\mathcal{A}_2^{qnil}$. Since $a^2b=aba$, therefore using the representation of $a$ and $b$ we obtain$$b_2=0,a_1b_1=b_1a_1\text{ and }a_2^2b_4=a_2b_4a_2.$$ Again because of $ab^2=bab$, it follows that $a_2b_4^2=b_4a_2b_4$. Since $aa^{nsd}b\in \mathcal{A}^{nsd}$ then by Lemma \ref{l3}, $b_4\in\mathcal{A}_2^{nsd}$. Therefore, for the idempotent $p_1=b_4b_4^d,$ 
  $b_4$ have the representation $\begin{bmatrix}
        b_{41} &0\\
        0 & b_{42}\\
    \end{bmatrix}_{p_1},$ where  $b_{41}\in (p_1\mathcal{A}_2p_1)^{inv}$ and $b_{42}\in(1-p_1)\mathcal{A}_2(1-p_1)^{qnil},$ correspondingly for $p_1$, $a_2$ have the representation $\begin{bmatrix}
        a_{21} & a_{22}\\
        a_{23} & a_{24}\\
    \end{bmatrix}_{p_1}.$  Now $a_2b_4^2=b_4a_2b_4$ implies $a_{23}=0$, then $a_2=\begin{bmatrix}
        a_{21} & a_{22}\\
        0 & a_{24}\\
    \end{bmatrix}_{p_1},$ where $a_{24}b_{42}^2=b_{42}a_{24}b_{42}$ and $b_{41}a_{21}=a_{21}b_{41}$. Since $a_{24},b_{42}\in((1-p_1)\mathcal{A}_2(1-p_1))^{qnil}$ with $a_{24}b_{42}^2=b_{42}a_{24}b_{42}$ and $a_{24}^2b_{42}=a_{24}b_{42}a_{24}$ then $a_{24}+b_{42}\in ((1-p_1)\mathcal{A}_2(1-p_1))^{qnil}$ by Lemma \ref{mylemma}. Hence according to Lemma \ref{l3}, in order to prove $a_2+b_4\in \mathcal{A}_2^{nsd}$ it enough to proof $a_{21}+b_{41}\in (p_1\mathcal{A}_2p_1)^{nsd}.$ Again since $a_{21}+b_{41}=b_{41}(p_1+(b_{41})^{-1}a_{21})$, therefore by Lemma \ref{c3}, if $p_1+(b_{41})^{-1}a_{21}\in  p_1\mathcal{A}_2p_1^{nsd}$  then  $a_{21}+b_{41}\in (p_1\mathcal{A}_2p_1)^{nsd}.$ Now \begin{align*}
        (p_1+(b_{41})^{-1}a_{21})-(p_1+(b_{41})^{-1}a_{21})^{n+1}\\=\big(p_1-\sum_{i=1}^{n+1}\binom{n+1}{i}b_{41}^{1-i}a_{21}^{i-1}\big)(b_{41})^{-1}a_{21}\in (p_1\mathcal{A}_2p_1)^{qnil}, 
    \end{align*}
    since $a_{21}\in(p_1\mathcal{A}_2p_1)^{qnil}$ and $b_{41}a_{21}=a_{21}b_{41}.$ Therefore $a+b\in\mathcal{A}^{nsd}$ if and only if $a_1+b_1\in (\mathcal{A}_1)^{nsd}$ according to Lemma \ref{l3}. Moreover $a_1b_1=b_1a_1$ and $aa^{nsd}b\in\mathcal{A}^{nsd}$ therefore for $p_2=b_1b_1^{nsd}$, $b_1,\text{ }a_1$ has the representation
    \begin{align*}
        b_1=\begin{bmatrix}
            b_{11} & 0\\
            0 & b_{12}\\
        \end{bmatrix}_{p_2} \hspace{0.2cm} \text{and} \hspace{0.2cm} a_1=\begin{bmatrix}
            a_{11}&0\\
            0 & a_{12}\\
        \end{bmatrix}_{p_2},
    \end{align*}
    where $b_{12}\in((1-p_2)\mathcal{A}_1(1-p_2))^{qnil}$ and $a_{12}\in((1-p_2)\mathcal{A}_1(1-p_2))^{nsd}$. Then by Theorem \ref{theo}, $a_{12}+b_{12}\in (1-p_2)\mathcal{A}_1(1-p_2)^{nsd}$. Hence $a_1+b_1\in \mathcal{A}_1^{nsd}$ if and only if $a_{11}+b_{11}\in p_2\mathcal{A}_1p_2^{nsd},$ and each of condition \ref{it1}-\ref{it3} is equivalent to $a_{11}+b_{11}\in (p_2\mathcal{A}p_2)^{nsd}.$ Conversely, if $a+b\in \mathcal{A}^{nsd}$ then $ a_{11}+b_{11}\in(p_2\mathcal{A}p_2)^{nsd}$ therefore conditions \ref{it1}-\ref{it3} holds. 
    \end{proof}
    
\end{theorem}
\begin{remark}
    H Chen, M Sheibani \cite{chen2022g} proved the g-Drazin invertibility of $a+b$ is equivalent to condition \ref{11}-\ref{12} given in Theorem \ref{2022}, under the setting $aba=\lambda a^2b=\lambda'ba^2$ and $bab=\alpha b^2a=\alpha' ab^2$, where $a,b\in\mathcal{A}^d$ and $\lambda,\lambda',\alpha,\alpha'$ are scalars. Theorem \ref{2022},  extended the results of H Chen, M Sheibani \cite{chen2022g} for the case of g$n$s-Drazin inverse when $\lambda=\alpha=1$, i.e. $aba=a^2b$ and $bab=ab^2.$
\end{remark}
\section{Additive results for p$n$s-Drazin inverse}\label{sec3}
In this section, our focus revolves around the concept of p$n$s-Drazin invertibility within the context of a Banach algebra.
Like Section \ref{sec1}, we establish some necessary and  sufficient conditions for the p$n$s-Drazin invertibility of an element in a Banach algebra. Using these necessary and sufficient conditions, we provide some additive results of p$n$s-Drazin inverse. \\First, we state some known results for the elements of Jacobson radical in a Banach algebra.
\begin{lemma}\cite{zou2017pseudo}\label{lemma21}
    Let $a,b\in \mathcal{A}$. 
    \begin{enumerate}
        \item If $a\in\mathcal{J}(\mathcal{A})$ or $b\in \mathcal{J}(\mathcal{A})$, then $ab,ba\in\mathcal{J}(\mathcal{A});$
        \item If $a\in \mathcal{J}(\mathcal{A})$ and $b\in \mathcal{J}(\mathcal{A})$, then $a+b\in \mathcal{J}(\mathcal{A}).$
    \end{enumerate} 
\end{lemma}
\begin{lemma}\cite{zou2017pseudo}\label{plusminus}
    Let $a,b\in \sqrt{\mathcal{J}(\mathcal{A})}$ with $ab=0$ or $ab=ba$. Then $a\pm b\in \sqrt{\mathcal{J}(\mathcal{A})}.$

\end{lemma}

The following lemma proves that the statement of Lemma   \ref{plusminus} holds under a condition weaker than commutativity.
\begin{lemma}
    Let $a,b\in \sqrt{\mathcal{J}(\mathcal{A})}$ with $ab^2=bab$ and $a^2b=aba$. Then $a\pm b\in \sqrt{\mathcal{J}(\mathcal{A})}.$
    \begin{proof}
        Since $ a,b \in \sqrt{\mathcal{J}(\mathcal{A})}, $ then $a^{k_1},b^{k_2}\in\mathcal{J}(\mathcal{A})$ for some $k_1,k_2\in\mathbb{N}.$ Therefore for $n>\text{max}\{k_1,k_2\}$, each monomial of $(a+b)^{3n}$, are an element of Jacobson radical. Hence $(a+b)^{3n}\in \mathcal{J}(\mathcal{A}),$ by Lemma \ref{lemma21}. Thus $a+b\in\sqrt{\mathcal{J}(\mathcal{A})}.$
    \end{proof}
\end{lemma}
\begin{proposition}\label{pC1}
    Let $n\in\mathbb{N}$ and $a\in\mathcal{A}$. Then $a^n-a^{2n}\in\sqrt{\mathcal{J}(\mathcal{A})}$ if and only if $a-a^{n+1}\in\sqrt{\mathcal{J}(\mathcal{A})}$.
    \begin{proof}
        If $a^n-a^{2n}\in\sqrt{\mathcal{J}(\mathcal{A})}$, then we have  $(a^n-a^{2n})^{k}\in J(\mathcal{A}),$ for some $k\in\mathbb{N}$. Since $$(a-a^{n+1})^{nk}=(a^n-a^{2n})^k(1-a^n)^{(n-1)k},$$ therefore $(a-a^{n+1})^{nk}\in \mathcal{J}(\mathcal{A})$ by Lemma \ref{lemma21}, and hence we obtain  $(a-a^{n+1})\in \sqrt{\mathcal{J}(\mathcal{A})}.$ \\Conversely, if $(a-a^{n+1})^k\in \mathcal{J}(\mathcal{A})$ for some $k\in \mathbb{N}$, then $$(a^n-a^{2n})^k=a^{(n-1)k}(a-a^{n+1})^k\in\mathcal{J}(\mathcal{A}).$$ Thus $a^n-a^{2n}\in \sqrt{\mathcal{J}(\mathcal{A})}.$
    \end{proof}
\end{proposition}
Now, we establish the main Theorem of this section.
\begin{theorem}\label{pT1}
    If $a\in\mathcal{A}$ and $n\in \mathbb{N}$, then $a$ is p$n$s-Drazin invertible if and only if  $a^n-a^{2n}\in\sqrt{\mathcal{J}(\mathcal{A})}.$
   
    \begin{proof}
        Suppose $a\in\mathcal{A}^{pnsD}$, then $(a^n-aa^{pnsD})^k\in\mathcal{J}(\mathcal{A}), \text{ for some }k\in\mathbb{N}.$ Now \begin{align*}
            a^n-a^{2n}&=a^n-aa^{pnsD}+aa^{pnsD}-a^{2n}\\&=(a^n-aa^{pnsD})+(aa^{pnsD}+a^n)(aa^{pnsD}-a^n).
        \end{align*}
        Therefore using Lemma \ref{plusminus}, we get $a^n-a^{2n}\in\sqrt{\mathcal{J}(\mathcal{A})}.$
        Conversely, if $a^n-a^{2n}\in \sqrt{\mathcal{J}(\mathcal{A})}$, then $a^n-a^{2n}\in\mathcal{A}^{qnil}$, thus $a$ is g-Drazin invertible and $a(1-aa^d)\in \mathcal{A}^{qnil}$. Furthermore $(a^n-1+aa^{d})((a^{d})^n-1+aa^{d})=1-(a^n-a^{n+1}a^{d})\in\mathcal{A}^{inv},$ hence $(a^n-1+aa^{d})\in\mathcal{A}^{inv}$. Moreover $$(a^n-aa^{d})(a^n-1+aa^{d})=a^{2n}-a^n\in \sqrt{\mathcal{J}(\mathcal{A})}.$$ Therefore $a^n-aa^{d}\in \sqrt{\mathcal{J}(\mathcal{A})}$, hence we have $a\in \mathcal{A}^{pnsD}.$
    \end{proof}

\end{theorem}
\begin{corollary}\label{pC2}
    Let $n\in\mathbb{N}$ and $a\in\mathcal{A}$. Then $a$ is p$n$s-Drazin invertible if and only if
     $a-a^{n+1}\in  \sqrt{\mathcal{J}(\mathcal{A})}.$
   
    \begin{proof}
        Follows from Theorem \ref{pT1} and Proposition \ref{pC1}. 
    \end{proof}
\end{corollary}
\begin{corollary}
    Let $n\in \mathbb{N}$. If $a\in\mathcal{A}^{pnsD}$ then $a\in\mathcal{A}^{p2nsD}.$ 
    \begin{proof}
        Since $a\in \mathcal{A}^{pnsD}$, therefore $a^n-a^{2n}\in\sqrt{\mathcal{J}(\mathcal{A})}.$ Furthermore $$a^{2n}-a^{4n}=(a^n-a^{2n})(a+a^{2n})\in \sqrt{\mathcal{J}(\mathcal{A})},$$ then by Theorem \ref{pT1}, we have $a\in\mathcal{A}^{p2nsD}.$
    \end{proof}
\end{corollary}
\begin{theorem}
    Let $a,b\in\mathcal{A}^{pnsD}$, for some $n\in\mathbb{N}$ such that $ab^2=bab$, then $ab\in \mathcal{A}^{pnsD}.$
    \begin{proof}
        It follows from equation \ref{eqn1} and Theorem \ref{pT1}.
    \end{proof}
\end{theorem}
Like in the case of g$n$s-Drazin inverse, now we investigate some additive properties for p$n$s-Drazin inverse. 
\begin{theorem}
    Let $a,b\in\mathcal{A}^{pnsD}$, for some $n\in\mathbb{N}$ such that $ab=0$, then $a+b\in \mathcal{A}^{pnsD}.$
    \begin{proof}
        Here $a,b$ are p$n$s-Drazin invertible, therefore $a-a^{n+1},b-b^{n+1}\in \sqrt{\mathcal{J}(\mathcal{A})}$. Since $ab=0$ then \begin{align*}
            (a+b)-(a+b)^{n+1}&=(a-a^{n+1})+(\sum\limits_{i=1}^{n}b^ia^{n+1-i})+(b-b^{n+1})\in\sqrt{\mathcal{J}(\mathcal{A})}.
        \end{align*} Therefore, in light of the Corollary \ref{pC2}, $a+b\in \mathcal{A}^{pnsD}$.
    \end{proof}
\end{theorem}

As Theorem \ref{theo}, we  establish the following result.  

\begin{theorem}
    Let $a,b\in\mathcal{A}^{pnsD}$, for some $n\in\mathbb{N}$ such that $ab=ba$, if $a+b\in \mathcal{A}^{pnsD}$ then $ab\left(\displaystyle\sum_{k=1}^n\binom{n+1}{k}a^{k-1}b^{n-k}\right)\in\sqrt{\mathcal{J}(\mathcal{A})}$ and if $ab\in \sqrt{\mathcal{J}(\mathcal{A})}$ then $a+b\in {(\mathcal{A})}^{pnsD}.$
\end{theorem}
As Lemma \ref{l3}, we get the next result for p$n$s-Drazin invertibility. 
\begin{lemma}\label{pT2}
    Let $x,y\in \mathcal{A}$ such that $x=\begin{bmatrix}
        a & c\\
        0 & b\\
    \end{bmatrix}_p$ and $y=\begin{bmatrix}
        b & 0\\
        c & a\\
    \end{bmatrix}_{(1-p)}$ where $p$ is an idempotent. Then for $n\in\mathbb{N}$
\begin{enumerate}
    \item if $a\in \mathcal{A}_1^{pnsD}$ and $b\in \mathcal{A}_2^{pnsD}$ then $x,y \in \mathcal{A}^{pnsD};$
    \item if $x\in \mathcal{A}^{pnsD}\left(y\in\mathcal{A}^{pnsD}\right)$ then $a\in\mathcal{A}_1^{pnsD}$ and $b\in \mathcal{A}_2^{pnsD}.$
\end{enumerate}

\begin{proof}
    \begin{enumerate}
        \item  Let $k_1$ and $k_2$ are two  positive integer satisfying $(a^n-a^{2n})^{k_1}\in \mathcal{J}(\mathcal{A})$ and $(b^n-b^{2n})^{k_2}\in\mathcal{J}(\mathcal{A})$, respectively. Choose $k=\text{max}\{k_1,k_2\}$, then \begin{equation*}
            (x^n-x^{2n})^{2k}\in\mathcal{J}(\mathcal{A}).
        \end{equation*} Therefore $x\in\mathcal{A}^{pnsD}.$ Similarly, it is true for $y$ also.
        \item  Since $x^n -x^{2n}\in\sqrt{\mathcal{J}(\mathcal{A})}$ therefore we have $b^n-b^{2n}\in\sqrt{\mathcal{J}(\mathcal{A}_2)}$  [\cite{zou2017pseudo}, Lemma 2.6]. Hence we get $b\in \mathcal{A}_2^{pnsD},$ by Theorem \ref{pT1}. Similarly we can prove $a\in\mathcal{A}_2^{pnsD}.$
    \end{enumerate}
\end{proof}
\end{lemma}
\begin{lemma}\label{pL1}
    Let $n\in\mathbb{N}$, $b\in \mathcal{A}^{pnsD}$ and  $a\in\sqrt{\mathcal{J}(\mathcal{A})}$ such that $$a(b^{\pi}-1)=0\hspace{0.2cm}\text{and}\hspace{0.2cm} b^{\pi}ab=0$$ then $a+b\in\mathcal{A}^{pnsD}.$
    \begin{proof}
        Since $b\in\mathcal{A}^{pnsD}$ therefore $b\in\mathcal{A}^{pD}$, then for the idempotent $p=bb^{pD}$, $b$ and $a$ has the representation  [\cite{zou2017pseudo}, Theorem 3.1] $$\begin{bmatrix}
           b_1 & 0\\
           0 & b_2\\
        \end{bmatrix}_p\text{ and }
        \begin{bmatrix}
            a_3 & a_1\\
            a_4 & a_2\\
        \end{bmatrix}_p,
        $$ where $b_1\in\mathcal{A}_1^{inv}$ and $b_2\in\sqrt{\mathcal{J}(\mathcal{A}_2)}$.  Now $a(b^\pi -1)=0 $ implies $ a=\begin{bmatrix}
            0 & a_1\\
            0 & a_2\\
        \end{bmatrix}_p$, where $a_2\in \sqrt{\mathcal{J}(\mathcal{A}_2)}$. Again \begin{equation*}
            b^\pi ab=0\implies a_2b_2=0,
        \end{equation*} and $a_2$, $b_2$ $\in\sqrt{\mathcal{J}(\mathcal{A}_2)}$, therefore $a_2+b_2\in\sqrt{\mathcal{J}(\mathcal{A}_2)}$. Since $b\in \mathcal{A}^{pnsD}$ then by Lemma \ref{pT2}, $b_1\in\mathcal{A}_2^{pnsD}$ hence $$a+b=\begin{bmatrix}
            b_1& a_1\\
            0 & a_2+b_2\\
        \end{bmatrix}_p,$$ is pns-Drazin invertible in light of Lemma \ref{pT2}.
    \end{proof}
\end{lemma}
\begin{theorem}
    Let $n\in\mathbb{N}$ and $a,b \in\mathcal{A}^{pnsD}$ such that $a(1-b^\pi)=0\hspace{0.2cm},b(1-a^\pi)=0$ and $b^\pi ab a^\pi=0$ then $a+b\in\mathcal{A}^{pnsD}.$
    \begin{proof}
      For the idempotent $p=aa^{pD}$, $a$ and $b$ has the representation  $$a=\begin{bmatrix}
          a_1 & 0\\
          0 & a_2\\
      \end{bmatrix}_p\text{ and } 
      \begin{bmatrix}
          b_3 & b_1\\
          b_4& b_2\\
      \end{bmatrix}_p,
      $$ where $a_1\in \mathcal{A}_1^{inv}$ and $a_2\in \sqrt{\mathcal{J}(\mathcal{A}_2)}$. Now $$b(1-a^\pi)=0\implies b=\begin{bmatrix}
          0 & b_1\\
          0& b_2\\
      \end{bmatrix}_p$$ and $$a(1-b^\pi)=0\implies a_2(1-b_2^{\pi})=0.$$ Again $$b^\pi aba^\pi=0\implies b_2^{\pi}a_2b_2=0.$$ Therefore by Lemma \ref{pL1} and Lemma \ref{pT2}, we have $a+b\in \mathcal{A}^{pnsD}.$
    \end{proof}
\end{theorem}
\begin{lemma}\label{psDl1}
   Let $n\in\mathbb{N}$. If $b\in\mathcal{A}^{pnsD}$ and $a\in\sqrt{\mathcal{J}(\mathcal{A})}$ such that $a(1-b^\pi)=0$ and $b^\pi ab=b^\pi ba$ then $a+b\in\mathcal{A}^{pnsD}.$
    \begin{proof}
        For the idempotent $p=bb^{pD}$, $b$ and $a$ has the representation   $$\begin{bmatrix}
           b_1 & 0\\
           0 & b_2\\
        \end{bmatrix}_p\text{ and }
        \begin{bmatrix}
            a_3 & a_1\\
            a_4 & a_2\\
        \end{bmatrix}_p,
        $$ where $b_1\in\mathcal{A}_1^{inv}$, $b_2\in\sqrt{\mathcal{J}(\mathcal{A}_2)}$.  Now $$a(1-b^\pi)=0\implies a=\begin{bmatrix}
            0 & a_1\\
            0 & a_2\\ 
        \end{bmatrix}_p \hspace{0.2cm}\text{and}\hspace{0.2cm}b^\pi ba=b^\pi ab\implies b_2a_2=a_2b_2.$$ Since $a_2,b_2\in \sqrt{\mathcal{J}(\mathcal{A})}$ and $a_2b_2=b_2a_2$, then by Lemma \ref{plusminus}, we get $a_2+b_2\in \sqrt{\mathcal{J}(\mathcal{A})}$. Therefore $a+b\in \mathcal{A}^{pnsD}.$
    \end{proof}
\end{lemma}
Similar to Theorem \ref{t3}, using Lemma \ref{psDl1}, we get the following result for p$n$s-Drazin invertibility.
\begin{theorem}
   Let $n\in\mathbb{N}$. If $a,b\in\mathcal{A}^{pnsD}$ such that $a(1-b^\pi)=0$, $b^\pi ba^\pi=b^\pi b$ and $b^\pi a^\pi ba=b^\pi a^\pi ab$ then $a+b\in \mathcal{A}^{pnsD}.$
\end{theorem}
Using Lemma \ref{pT2}, the next result follows similarly as Theorem \ref{Theor}. 

\begin{theorem}
    Let $n\in\mathbb{N}$. If $a,b$ are two p$n$s-Drazin invertible elements in $\mathcal{A}$ such that $aba^{\pi}=0$.  Then $a+b\in \mathcal{A}^{pnsD}$ if and only if $aa^{pnsD}(a+b)\in \mathcal{A}^{pnsD}.$
\end{theorem}
As Theorem \ref{2022}, using the properties of Jacobson radical in place of quasinilpotent elements, we get the next result.
\begin{theorem}
    Let $a,b\in \mathcal{A}^{pnsD}$ for some $n\in \mathbb{N}$, such that $a^2b=aba,ab^2=bab$. Then  the followings are equivalent: 
    \begin{enumerate}
        \item  $a+b\in \mathcal{A}^{pnsD}$;
        \item \label{it11} $aa^{pnsD}(a+b)bb^{pnsD}\in \mathcal{A}^{pnsD}$;
        \item\label{it22} $aa^{pnsD}(a+b)\in \mathcal{A}^{pnsD};$
        \item\label{it33} $(a+b)bb^{pnsD}\in \mathcal{A}^{pnsD}.$
        
    \end{enumerate}
\end{theorem}
\section{Additive properties for weighted strong Drazin inverse}\label{sec4}
\indent For a non-zero element $w\in\mathcal{A}$, $w$-product on $\mathcal{A}$ is defined as $a*b=awb$, for $a,b\in\mathcal{A}$. With this $w$-product and the $w$-norm $\|a\|_w=\|a\|\|w\|$, $\mathcal{A}$
becomes complex Banach algebra, which will be denoted by $\mathcal{A}_w$. By $r_w(a)\text{ and }r(a)$ we will denote the spectrum of $a$ in $\mathcal{A}_w$ and $\mathcal{A}$, respectively. For $w\in\mathcal{A}-\{0\}$, Dijana Mosić\cite{mosic2021generalized} defined weighted generalized strongly Drazin inverse(wgs-Drazin inverse) and weighted pseudo strongly Drazin inverse(wps-Drazin inverse) as follows, an element $a\in\mathcal{A}$ is wgs-Drazin invertible if there exists $x\in\mathcal{A}$ such that $$x*a=a*x,\text{ }x*a*x=x\text{ and }(a-a*x)\in \mathcal{A}_w^{qnil}.$$

Similarly, an element $a\in\mathcal{A}$ is wps-Drazin invertible if there exists $x\in\mathcal{A}$ such that $$x*a=a*x,\text{ }x*a*x=x\text{ and }(a-a*x)^{*k}\in \mathcal{J}(\mathcal{A})\text{ for some }k\in\mathbb{N}.$$ 
Like g$n$s and p$n$s-Drazin inverse for $n\in\mathbb{N},$ now we define weighted generalized n-strong Drazin inverse(wg$n$s-Drazin inverse) and weighted pseudo n-strong Drazin inverse(wp$n$s-Drazin inverse), for some non-zero weight $w\in\mathcal{A}.$
\begin{definition}

  Let $w\in\mathcal{A}-\{0\}\text{ and } n\in\mathbb{N}$:   
  \begin{enumerate}
      \item $a\in\mathcal{A}$ is wg$n$s-Drazin invertible, if there exists $x\in\mathcal{A}$ satisfying $$x*a=a*x,\text{ }x*a*x=x\text{ and }a^{*n}-a*x\in\mathcal{A}^{qnil}_w.$$ Here $x$ is wg$n$s-Drazin inverse of $a$, denoted by $a^{nsd,w}.$
      \item Similarly $a\in\mathcal{A}$ is wp$n$s-Drazin invertible, if there exists $y\in \mathcal{A}$ such that
      \begin{multline*}
        y*a=a*y,\text{ }y*a*y=y\text{ and } \left(a^{*n}-a*y\right)^{*k}\in\mathcal{J}(\mathcal{A})\text{ for some }k\in\mathbb{N}.  
      \end{multline*}
      Here $y$ is wp$n$s-Drazin inverse of $a$, denoted by $a^{pnsD,W}.$
  \end{enumerate}
\end{definition}
The collection of wg$n$s and wp$n$s-Drazin invertible elements in the Banach algebra $\mathcal{A}$ are symbolized by $\mathcal{A}^{nsd,w}$ and $\mathcal{A}^{pnsD,w}$, respectively. Notice that $a\in\mathcal{A}^{nsd,w}$ means $a$ is g$n$s-Drazin invertible as an element in the algebra $\mathcal{A}_w$; therefore wg$n$s-Drazin inverse is unique. Similarly wp$n$s-Drazin inverse is also unique. 
\\ \indent Next lemma is very much useful in the study of wg$n$s and wp$n$s-Drazin inverses.
\begin{lemma}\cite{dajic2007weighted}\label{spec}
    For $w,a\in\mathcal{A}$, $r(aw)=r_w(a)=r(wa)$. Therefore $a\in\mathcal{A}_w^{qnil}$ if and only if $aw\in\mathcal{A}^{qnil}$ if and only if $wa\in\mathcal{A}^{qnil}.$
\end{lemma}
Now we characterize the wg$n$s and wp$n$s-Drazin invertibility of an element in a Banach algebra in the following theorems.
\begin{theorem}\label{wgns}
   Let $w\in \mathcal{A}-\{0\}$ and $n\in\mathbb{N}.$ Then for $a\in\mathcal{A}$ the followings are equivalent:
    \begin{enumerate}
       \item $a\in\mathcal{A}^{nsd,w};$\label{1}
        \item $wa\in\mathcal{A}^{nsd};$\label{3}
       \item $aw\in\mathcal{A}^{nsd}.$\label{2}
     
      The wg$n$s-Drazin inverse will satisfies $a^{nsd,w}=\left((aw)^{nsd}\right)^2a$.
    \end{enumerate}
    \begin{proof}
        \ref{1}$\implies$\ref{2} Let $c=a^{nsd,w}$, then we have $a*c=c*a$ which implies $ awcw=cwaw$ and $c*a*c=c$ implies $ cwawcw=cw$. Moreover $$a^{*n}-a*c\in\mathcal{A}^{qnil}_w\implies (aw)^n-awcw\in \mathcal{A}^{qnil}\text{ by Lemma \ref{spec}}.$$ Hence we get $aw\in\mathcal{A}^{nsd}.$
       \\ \indent \ref{2}$\implies$\ref{1} Put $c=\left((aw)^{nsd}\right)^2a$, then by definition $a*c=c*a\text{ and }c*a*c=c$. Since $(a^{*n}-a*c)w=(aw)^n-(aw)(aw)^{nsd}$, therefore by Lemma \ref{spec}, $a^{*n}-c*a\in\mathcal{A}^{qnil}_w.$ Hence we have $a\in\mathcal{A}^{nsd,w}.$  
        \\ \indent Similarly \ref{1} and \ref{3} is also equivalent.
        
        \end{proof}
\end{theorem}

\begin{theorem}\label{wgns1}
    If $w\in\mathcal{A}-\{0\}$ and $n\in\mathbb{N}$, then $a\in\mathcal{A}$ is wg$n$s-Drazin invertible if and only if $a^{*n}-a^{*2n}\in\mathcal{A}_w^{qnil}.$
    \begin{proof}
        If $a\in\mathcal{A}^{nsd,w}$ then by Theorem \ref{wgns}, $aw\in\mathcal{A}^{nsd}.$ Now
        \begin{align*}
            (a^{*n}-a^{*2n})w&=(aw)^n-(aw)^{2n}\\&=(aw)^n-aw(aw)^{nsd}+\left(aw(aw)^{nsd}\right)^2-(aw)^{2n}
            \\&=\left((aw)^n-aw(aw)^{nsd}\right)+\left((aw)^n-aw(aw)^{nsd}\right)\left((aw)^n+aw(aw)^{nsd})\right).
        \end{align*}
        Therefore using Lemma \ref{mylemma}, we have $(a^{*n}-a^{*2n})w\in\mathcal{A}^{qnil}$ and by Lemma \ref{spec}, $a^{*n}-a^{*2n}\in\mathcal{A}_w^{qnil}.$
        Conversely, if $a^{*n}-a^{*2n}\in\mathcal{A}_w^{qnil}$ then using Lemma \ref{spec}, $(a^{*n}-a^{*2n})w=(aw)^n-(aw)^{2n}\in\mathcal{A}^{qnil},$ then  $aw\in\mathcal{A}^{nsd}$ in light of Theorem \ref{t1}. Hence we get $a\in\mathcal{A}^{nsd,w}$ from Theorem \ref{wgns}. 
    \end{proof}
\end{theorem}
As Theorem \ref{wgns}, we get the following theorem.

\begin{theorem}\label{wpnsD}
   Let $w\in \mathcal{A}-\{0\}$ and $n\in\mathbb{N}$. Then for $a\in\mathcal{A}$ the following statements are equivalent:
    \begin{enumerate}
       \item $a\in\mathcal{A}^{pnsD,w};$
       \item $wa\in\mathcal{A}^{pnsD};$
       \item $aw\in\mathcal{A}^{pnsD}.$
      
      The wp$n$s-Drazin inverse will satisfies $a^{pnsD,w}=\left((aw)^{pnsD}\right)^2a$.
    \end{enumerate}
\end{theorem}
Similarly, as Theorem \ref{wgns1}, we get the following characterization of wp$n$s-Drazin inverse.

\begin{theorem}
    If $w\in\mathcal{A}-\{0\}$ and $n\in\mathbb{N}$. Then $a\in\mathcal{A}$ is wp$n$s-Drazin invertible if and only if $(a^{*n}-a^{*2n})^{*k}\in\mathcal{J}(\mathcal{A})$ for some $k\in\mathbb{N}.$
    \begin{proof}
        If $a\in\mathcal{A}^{pnsD,w}$ then by Theorem \ref{wpnsD}, $aw\in\mathcal{A}^{pnsD}.$ Now
        \begin{align*}
            (a^{*n}-a^{*2n})w&=(aw)^n-(aw)^{2n}\\&=(aw)^n-aw(aw)^{pnsD}+\left(aw(aw)^{pnsD}\right)^{2n}-(aw)^{2n}
            \\&=\left(aw-aw(aw)^{pnsD}\right)+\\&\left((aw)^n-(aw(aw)^{pnsD})^n\right)\left((aw)^n+(aw(aw)^{pnsD})^{2n})\right).
        \end{align*}
        Then using Lemma \ref{plusminus}, $(a^{*n}-a^{*2n})w\in\sqrt{\mathcal{J}(\mathcal{A})}$, hence $(a^{*n}-a^{*2n})^{*k}\in\mathcal{J}(\mathcal{A})$ for some $k\in\mathbb{N}.$
        On the other side if $(a^{*n}-a^{*2n})^{*k}\in\mathcal{J}(\mathcal{A})$ then $(a^{*n}-a^{*2n})^{*k}w=\left((aw)^n-(aw)^{2n}\right)^k\in\sqrt{\mathcal{J}(\mathcal{A})},$ then by Theorem \ref{pT1}, $aw\in\mathcal{A}^{pnsD}$. Hence in light of the Theorem 
 \ref{wpnsD}, $a\in\mathcal{A}^{pnsD,w}.$ 
    \end{proof}
\end{theorem}
Since $x\in\mathcal{A}^{nsd,w}\left(\mathcal{A}^{pnsD,w}\right)$ implies $x$ is g$n$s(p$n$s)-Drazin invertible in the Banach algebra $\mathcal{A}_w$, therefore all the additive results developed in Section \ref{sec1}(Section \ref{sec3}), is true for the sum of two elements in the algebra of $\mathcal{A}_w$ also. Hence all those results follow in the case of wg$n$s(wp$n$s)-Drazin inverse also.

\subsection*{Acknowledgment}

\end{document}